\documentclass[12pt]{article}
\begin{document}
\newpage
\pagestyle{empty}
\setcounter{page}{0}
%
\newfont{\twelvemsb}{msbm10 scaled\magstep1}
\newfont{\eightmsb}{msbm8}
\newfont{\sixmsb}{msbm6}
\newfam\msbfam
\textfont\msbfam=\twelvemsb
\scriptfont\msbfam=\eightmsb
\scriptscriptfont\msbfam=\sixmsb
\catcode`\@=11
\def\Bbb{\ifmmode\let\next\Bbb@\else
  \def\next{\errmessage{Use \string\Bbb\space only in math mode}}\fi\next}
\def\Bbb@#1{{\Bbb@@{#1}}}
\def\Bbb@@#1{\fam\msbfam#1}
\newfont{\twelvegoth}{eufm10 scaled\magstep1}
\newfont{\tengoth}{eufm10}
\newfont{\eightgoth}{eufm8}
\newfont{\sixgoth}{eufm6}
\newfam\gothfam
\textfont\gothfam=\twelvegoth
\scriptfont\gothfam=\eightgoth
\scriptscriptfont\gothfam=\sixgoth
\def\frak{\frak@}
\def\frak@#1{{\fam\gothfam{{#1}}}}
\def\frak@@#1{\fam\gothfam#1}
\catcode`@=12
%
%
%
\def\CC{{\Bbb C}}
\def\NN{{\Bbb N}}
\def\QQ{{\Bbb Q}}
\def\RR{{\Bbb R}}
\def\ZZ{{\Bbb Z}}
\def\cA{{\cal A}}          \def\cB{{\cal B}}          \def\cC{{\cal C}}
\def\cD{{\cal D}}          \def\cE{{\cal E}}          \def\cF{{\cal F}}
\def\cG{{\cal G}}          \def\cH{{\cal H}}          \def\cI{{\cal I}}
\def\cJ{{\cal J}}          \def\cK{{\cal K}}          \def\cL{{\cal L}} 
\def\cM{{\cal M}}          \def\cN{{\cal N}}          \def\cO{{\cal O}}
\def\cP{{\cal P}}          \def\cQ{{\cal Q}}          \def\cR{{\cal R}} 
\def\cS{{\cal S}}          \def\cT{{\cal T}}          \def\cU{{\cal U}}
\def\cV{{\cal V}}          \def\cW{{\cal W}}          \def\cX{{\cal X}}
\def\cY{{\cal Y}}          \def\cZ{{\cal Z}}
\def\qed{\hfill \rule{5pt}{5pt}}
\def\id{\mbox{id}}
\def\ggo{{\frak g}_{\bar 0}}
\def\uqggo{\cU_q({\frak g}_{\bar 0})}
\def\uqggp{\cU_q({\frak g}_+)}
\def\typeA{{\em type $\cA$}}
\def\typeB{{\em type $\cB$}}
\newtheorem{lemma}{Lemma}
\newtheorem{prop}{Proposition}
\newtheorem{theo}{Theorem}
%
%
\newcommand{\norm}[1]{{\protect\normalsize{#1}}}
\newcommand{\LAP}
{{\small E}\norm{N}{\large S}{\Large L}{\large A}\norm{P}{\small P}}
\newcommand{\sLAP}{{\scriptsize E}{\footnotesize{N}}{\small S}{\norm L}$
${\small A}{\footnotesize{P}}{\scriptsize P}}
\def\logolapin{
  \raisebox{-1.2cm}{\epsfbox{/lapphp8/keklapp/ragoucy/paper/enslapp.ps}}}
\def\logolight{{\bf{{\large E}{\Large N}{\LARGE S}{\huge L}{\LARGE
        A}{\Large P}{\large P}} }}
\def\logoenslapp{\logolight}
%
%
%
\hbox to \hsize{
\hss
\begin{minipage}{5.2cm}
  \begin{center}
    {\bf Groupe d'Annecy\\ \ \\
      Laboratoire d'Annecy-le-Vieux de Physique des Particules}
  \end{center}
\end{minipage}
\hfill
\logoenslapp
\hfill
\begin{minipage}{4.2cm}
  \begin{center}
    {\bf Groupe de Lyon\\ \ \\
      {\'E}cole Normale Sup{\'e}rieure de Lyon}
  \end{center}
\end{minipage}
\hss}

\vspace {.3cm}
\centerline{\rule{12cm}{.42mm}}
\vfill
\vfill
\begin{center}

  {\LARGE {\bf {\sf Centre and Representations of ${\cal
          U}_{q}(sl(2|1))$ at Roots of Unity }}} \\[1cm]

\vfill

{\large B. Abdesselam$^{\dagger,}$
\footnote{abdess@orphee.polytechnique.fr,

  ~$^2$arnaudon@lapp.in2p3.fr

  ~$^3$bauer@spht.saclay.cea.fr

  ~$^4$Permanent address: Service de Physique Th{\'e}orique,
  C.E.A. Saclay, F-91191, Gif-sur-Yvette, France.

  ~$^5$Laboratoire Propre du CNRS UPR A.0014.

  ~$^6$URA 14-36 du CNRS, associ{\'e}e {\`a} l'E.N.S. de Lyon et {\`a}
  l'Universit{\'e} de Savoie.

  ~$^7$Partially supported by European Community Contract
  ERBCHRXCT920069.}
  D. Arnaudon $^{\ddagger,2}$
  and M. Bauer$^{\ddagger,3,4,7}$}

\vfill

{\em $^{\dagger}$
  Centre de Physique Th{\'e}orique$^{^5}$, 
  {\'E}cole Polytechnique 
  91128 Palaiseau Cedex,
  France.

$^\ddagger$ \LAP$^6$, Chemin de Bellevue BP 110,
74941 Annecy-le-Vieux Cedex, France.}

\end{center}

\vfill

\begin{abstract}

Quantum groups at roots of unity have the property that their centre
is enlarged. Polynomial equations relate the standard deformed
Casimir operators and the new central elements. 
These relations are important from a physical point of view since they
correspond to relations among quantum expectation values of
observables that have to be satisfied on all physical states. 
In this paper, we
establish these relations in the case of the quantum Lie
superalgebra ${\cal U}_{q}(sl(2|1))$. In the course of the argument, we
find and use a set of representations such that any relation
satisfied on all the representations of the set is true in 
${\cal U}_{q}(sl(2|1))$. This set is a subset of the set of all the
finite dimensional   irreducible representations of 
${\cal U}_{q}(sl(2|1))$, that we classify and describe explicitly.

\end{abstract}

\vfill
\vfill

\rightline{\LAP-A-583/96}
\rightline{RR 443.0396}
\rightline{T96/019}

\newpage
\pagestyle{plain}

\section{Introduction\label{sect:introduction}}

Classical and quantum Lie superalgebras and their representations play 
respectively an important role in the understanding and exploitation of the 
classical and $q$-deformed supersymmetry in physical systems. A complete 
classification of the finite-dimensional simple classical Lie superalgebras 
over $\CC$ has been given by Kac \cite{Kaci,Kacii} and Scheunert
\cite{ScheuSuper}. The  
corresponding irreducible representations fall in two series called
typical and atypical.

Irreducible representations of the quantum analogue of superalgebras are 
studied intensively when $q$ is not root of 
unity in \cite{PSVdJ,PalTol,PalStoa,Zhanga}.  

A complete classification of finite dimensional 
irreducible representations of unrestricted quantum algebras  
for $q$ a root of unity exits only
in the particular case of ${\cal U}_{q}(sl(2))$ \cite{RA}.
Partial classifications exist also in the case of 
${\cal U}_{q}(sl(3))$, in \cite{DobStAndrews} for 
the restricted case and in \cite{Asutroisq} for periodic representations.
Much progress towards a complete classification in the general case of
$\cU_q(\cG)$ for $\cG$ a simple Lie algebra was done in \cite{DK,DKP}.

The classification of finite dimensional 
irreducible representations of $\cU_q(osp(1|2))$ for any $q$
parallels the 
${\cal U}_{q}(sl(2))$ case \cite{GCX,Kobayashi,PalStob,Zhangb}.
The only other fully understood case is ${\cal U}_{q}(sl(2|1))$
\cite{Zhanga,Zhangc}. 

Our main goal in this paper is the structure of the centre of ${\cal
  U}_{q}(sl(2|1))$ when $q$ is a root of unity. Complete sets of
representations - to be defined below - give a convenient way to prove
relations in this centre. Their construction involves a detailed
knowledge of matrix elements of the finite dimension irreducible
representations, whose classification is given below with emphasis on
what is needed for the rest of the paper. 

In Section \ref{sect:definitions}, 
we give the definition of the quantum superalgebra 
${\cal U}_{q}(sl(2|1))$ and the expression of  central elements. 
Generalities on the finite dimensional irreducible 
representations of ${\cal U}_{q}(sl(2|1))$ are presented in section
\ref{sect:generalities}. 
In section \ref{sect:sl2}, we recall some useful results on
$\cU_q(gl(2))$ at roots of unity and we give complete sets of
irreducible representations for this quantum algebra: expressions in
the universal quantum enveloping 
algebra that vanish on such sets, vanish identically.
In section \ref{sect:classification}, we classify the finite
dimensional irreducible representations of ${\cal U}_{q}(sl(2|1))$. 
In section \ref{sect:complete_set}, we present complete sets of
representations corresponding to infinite subsets of the set of
continuous parameters. All the representations of these complete sets
have the same dimension, unlike in the classical case \cite{ABP}.
Finally, in section \ref{sect:proof_rel}, we prove the relations in
the centre using our complete set of irreducible representations. 

\smallskip

\section{Quantum superalgebra ${\cal U}_{q}(sl(2|1))$ and its centre
  \label{sect:definitions}}

The superalgebra ${\cal U}_{q}(sl(2|1))$ is the associative 
superalgebra over $\CC$ with generators $k_{1}=q^{h_{1}}$, 
$k_{1}^{-1}=q^{-h_{1}}$, $k_{2}=q^{h_{2}}$, $k_{2}^{-1}=q^{-h_{2}}$, 
$e_{1}$, $e_{2}$, $f_{1}$, $f_{2}$ and  relations
\begin{eqnarray}
  && k_1 k_2 = k_2 k_1 \;, \\
  && k_i e_j k_i^{-1} = q^{a_{ij}} e_j\;,   \qquad\qquad
     k_i f_j k_i^{-1} = q^{-a_{ij}} f_j\;,   \label{eq:kifj}  \\
  && e_1 f_1 - f_1 e_1 = \frac{k_1-k_1^{-1}}{q-q^{-1}}\;,   \qquad\qquad 
     e_2 f_2 + f_2 e_2 = \frac{k_2-k_2^{-1}}{q-q^{-1}}\;,  \label{eq:e2f2}  \\
  && [e_1,f_2] = 0\;,  \qquad\qquad
     [e_2,f_1] = 0\;,  \label{eq:e2f1} \\
  && e_2^2 = f_2^2 = 0  \;, \label{eq:f2}  \\
  && e_1^2 e_2 - (q+q^{-1}) e_1 e_2 e_1+ e_2 e_1^2 = 0  \;,
  \label{eq:serre1}  \\ 
  && f_1^2 f_2 - (q+q^{-1}) f_1 f_2 f_1+ f_2 f_1^2 = 0  \;.
  \label{eq:serre2}   
\end{eqnarray}
The two last equations are called the Serre relations.
The matrix $(a_{ij})$ is the distinguished Cartan matrix of $sl(2|1)$,
i.e.  
\begin{equation}
  (a_{ij}) = 
  \left( 
    \begin{array}{cc} 2 & -1 \\ 
                     -1 & 0 
    \end{array}
  \right)
  \label{eq:Cart_mat}
\end{equation}

The $\ZZ_{2}$-grading in ${\cal U}_{q}(sl(2|1))$ is uniquely defined by the 
requirement that the only odd generators are $e_{2}$ and $f_{2}$, i.e.
\begin{eqnarray}
  &&\deg\;(k_{1})=\deg\;(k_{2})=0 \;, \nonumber\\
  &&\deg\;(k_{1}^{-1})=\deg\;(k_{2}^{-1})=0 \;, \nonumber\\
  &&\deg\;(e_{1})=\deg(f_{1})=0 \;, \nonumber\\
  &&\deg\;(e_{2})=\deg(f_{2})=1 \;. 
  \label{eq:degree}
\end{eqnarray}
We will not use the (standard) co-algebra structure in the following. 

\medskip

Define
\begin{equation}
  e_{3}=e_{1}\;e_{2}-q^{-1}e_{2}\;e_{1} \qquad \mbox{and} \qquad
  f_{3}=f_{2}\;f_{1}-q\;f_{1}\;f_{2} \;.
  \label{eq:e3etf3}
\end{equation}
The quantum Serre relations become
\begin{eqnarray}
  &&e_{1}\;e_{3}=q\;e_{3}\;e_{1} \;,\nonumber\\
  &&f_{3}\;f_{1}=q^{-1}f_{1}\;f_{3} \;. 
  \label{eq:serre3}
\end{eqnarray}
Furthermore,
\begin{eqnarray}
  && e_{2}\;e_{3}=-q\;e_{3}\;e_{2} \;, \nonumber\\
  && f_{3}\;f_{2}=-q^{-1}f_{2}\;f_{3} \;,
  \label{eq:e2f3}
\end{eqnarray}
and
\begin{eqnarray}
  && e_{3}\;f_{3}+f_{3}\;e_{3} =
  {k_{1}\;k_{2}-k_{1}^{-1}\;k_{2}^{-1} \over q-q^{-1}} \;, \nonumber\\
  && e_{3}^{2} = f_{3}^{2}=0 \;.
  \label{eq:e3f3}
\end{eqnarray}

\medskip

In the following, we will use the conventional notation 
\begin{equation}
{}[x] \equiv \frac{q^x-q^{-x}}{q-q^{-1}} \;.
  \label{eq:q_nombre}
\end{equation}

When $q$ is not a root of unity,
the centre of ${\cal U}_{q}(sl(2|1))$ is generated by the elements
$\cC_p$, $p\in\ZZ$, where
\begin{eqnarray}
  &{\cal C}_{p}=k_{1}^{2p-1}k_{2}^{4p-2}
  & \bigg. 
  (q-q^{-1})^2 
  \bigg\lbrace  [h_{1}+h_{2}+1][h_{2}]
  \nonumber\\
  && \bigg. \qquad  
  - f_{1}e_{1}+f_{2}e_{2}([h_{1}+h_{2}]q^{1-2p}-[h_{1}+h_{2}+1])
  \nonumber\\
  && \bigg. \qquad  +
  f_{3}e_{3}([h_{2}-2]q^{1-2p}-[h_{2}-1])
  \nonumber\\
  && \bigg. \qquad   +
  (q-q^{-1}) 
  q^{-1-p}[p]f_{3}e_{2}e_{1}k_{2} 
  + (q-q^{-1})  
  q^{2-p}f_{1}f_{2}e_{3}
  k_{2}^{-1}[p-1] \nonumber\\
  && \bigg. \qquad  +
  (q-q^{-1})^2 
  q^{1-2p}[p][p-1]f_{2}f_{3}e_{3}e_{2} \bigg\rbrace \;.
  \label{eq:cas_sl21}
\end{eqnarray}
They satisfy the relations 
\begin{equation}
  {\cal C}_{p_1} {\cal C}_{p_2}   ={\cal C}_{p_3} {\cal C}_{p_4} \qquad
  \hbox{if}
  \qquad p_1 + p_2 = p_3 + p_4 \;.
  \label{eq:rel1}
\end{equation}
The fact that the centre was not finitely generated in the classical
case was known since \cite{KacSuperRep,ScheuEig}. 
The explicit expression of a set of generators of the centre, together
with the relations, was given in
\cite{ABP} in the classical case and in \cite{ACF} in the quantum
case.

In this paper, we consider the case when $q$ is a root of unity. 
Let $l$ the smallest integer 
such that $q^{l}=1$. We define:
\begin{equation}
  l'=\cases{ l & if $l$ is odd  \cr
       &  \cr
       l/2 & if $l$ is even \cr}
  \label{eq:lprime}
\end{equation}
the elements $z_{i}\equiv k_{i}^{l}$, $x_{1}\equiv e_{1}^{l}$ and
$y_{1}\equiv f_{1}^{l}$ also belong to the centre. 

\begin{prop}
When $l$ is odd, the central elements
$z_1$, $z_2$, $x_1$, $y_1$ and $\cC_p$, $p\in\ZZ$
satisfy the relations
\begin{eqnarray}
  {\cal C}_{p+l} & = & z_1^2 z_2^4 {\cal C}_{p}, \nonumber\\
  {\cal C}_{p+1}^{l} & = & z_1^2 z_2^4 {\cal C}_{p}^{l},
  \nonumber\\ 
  \cP_l(\cC_1,\cdots,\cC_l) &\equiv &
  ({\cal C}_{1}+1)^{l}-1 + \sum_{m\geq 2 \;\; n\geq 0 \atop m+n \leq l}
  {\cal C}_{m}{\cal C}_{1}^{n}{l\over m-1}\biggl({m+n-1 \atop n+1}\biggr)
  \biggl({l-m \atop n}\biggr) \nonumber\\
  &=&
  \left( 1-z_1^2 z_2^2 \right)\left(z_2^2 -1\right) -
    (q-q^{-1})^{2l} 
    z_1^2 z_2^4 y_{1}x_{1} \;. 
  \label{eq:rel2}
\end{eqnarray}
\end{prop}
The first two relations follow from the expression of ${\cal C}_{p}$ and from
(\ref{eq:rel1}). 
The third relation will be proved using complete sets
of representations of
${\cal U}_{q}(sl(2|1))$. Furthermore, there is no other independent
polynomial relation. 

\smallskip

\section{Generalities on finite dimensional irreducible
  representations
  \label{sect:generalities}}

 Let us consider a finite dimensional irreducible left 
module $M$ over ${\cal U}_{q}(sl(2|1))$.

\begin{itemize}
\item The generators $k_{1}$ and $k_{2}$ are 
  simultaneously diagonalizable on the module $M$. 

\smallskip

\item Since $e_{2}^{2}=0$ and $\dim M < \infty$, there 
  exists a subspace $V \subset M$ annihilated by $e_{2}$, i.e. 
  \begin{equation}
    \forall\; v \in V, \qquad \qquad e_{2}\;v=0 \;.
  \end{equation}
  
  \smallskip
  
\item Since $e_{3}^{2}=0$ and $e_{2}\;e_{3}=-q\;e_{3}\;e_{2}$, 
  there exists $V_{0} \subset V$ annihilated by $e_{3}$, i.e.
  \begin{equation}
    \forall\; v \in V_{0}, \qquad \qquad
    e_{2}\;v=e_{3}\;v=0.   
  \end{equation}
  
  \smallskip

\item Because of (\ref{eq:kifj}) the subspace $V_{0}$ is stable 
  by left multiplication by $k_{1}$ and $k_{2}$.
  
  \smallskip
  
\item Because of (\ref{eq:e3etf3}) and (\ref{eq:serre3}) 
  the  subspace $V_{0}$ is stable 
  by left multiplication by $e_{1}$.
  
  \smallskip
  
\item Because of $e_{2}\;f_{1}=f_{1}\;e_{2}$ and $e_{3}\;f_{1}-
  f_{1}\;e_{3}=-e_{2}\;k_{1}^{-1}$ the  subspace $V_{0}$ is stable by left 
  multiplication by $f_{1}$.
  
\end{itemize}

Let $\ggo\simeq gl(2)$ be the even subalgebra of $sl(2|1)$.
The algebra 
${\cal U}_{q}(\ggo)$ is generated by $e_{1}$,
 $f_{1}$, $k_{1}$ and $k_{2}$.

The module $V_{0}$ is then a ${\cal U}_{q}(\ggo)$ submodule of $M$. 
It is simple (as a ${\cal U}_{q}(\ggo)$ module), since 
any submodule of $V_{0}$ would generate a proper submodule of $M$ by 
left action of ${\cal U}_{q}(sl(2|1))$. 
As a consequence of the simplicity of 
$V_{0}$, the element $k_{1}\;k_{2}^{2}$ (the $U(1)$ generator) is 
represented by a scalar on $V_{0}$.

\medskip

Let $\uqggp$ be the subalgebra of $\cU_q(sl(2|1))$ generated by 
$e_{1}$, $f_{1}$, $k_{i}$ and $e_2$. 
The subspace $V_0$ is also a $\uqggp$-module, annihilated by $e_2$.

\medskip

{}From $V_0$ considered as an $\uqggp$-module, one can construct an
induced $\cU_q(sl(2|1))$-module $M'=\cU_q(sl(2|1)) \otimes_{\uqggp}
V_0$. Then $M$ is  equal to $M'$ if $M'$ is simple, or to the quotient
of $M'$ by its maximal submodule otherwise. 

\medskip

Since we already know that each finite dimensional irreducible 
representation of $\cU_q(sl(2|1))$
is associated to one finite dimensional irreducible 
representation $V_0$ of ${\cal U}_{q}(\ggo)$, we will
construct the classification of the former in terms of the latter.  
As we will see, the correspondence is one-to-one. 
We now need some results on $\cU_q(gl(2))$ at roots of unity.

\section{$\cU_q(gl(2))$ at roots of unity \label{sect:sl2}}
 
\subsection{Centre of $\cU_q(gl(2))$}
The elements $k_1 k_2^2$, $k_1^{l'}$, $k_2^l$, $e_1^l$, $f_1^l$ 
$e_1^{l'}f_1^{l'}$, $f_1^{l'}e_1^{l'}$ 
are central in $\cU_q(gl(2))$. 
The $q$-deformed quadratic Casimir operator is
\begin{equation}
\cC_{\cU_q(gl(2))} = qk_1 + q^{-1}k_1^{-1} + (q-q^{-1})^{2} f_1 e_1 
  \label{eq:cas_gl2}\;.
\end{equation}
When $l$ is odd,
the centre of $\cU_q(gl(2))$ is actually the algebra defined by the
generators $k_1 k_2^2$, 
$k_1^{l}$, $e_1^l$, $f_1^l$, 
$\cC_{\cU_q(gl(2))}$ and the relation 
\cite{KerDipl}
\begin{equation}
  2
  P_{l}\left(\cC_{\cU_q(gl(2))}/2\right)  
  = k_1^{l} +  k_1^{-l} + (q-q^{-1})^{2l} f_1^{l} e_1^{l}  \;.
  \label{eq:rel_gl2}
\end{equation}
The polynomial 
$P_{l}$ is the first kind Chebychev Polynomial of degree $l$ defined by
\begin{equation}
  P_l(\cos x) = \cos (lx) \;.
  \label{eq:chebi}
\end{equation}

\subsection{Finite dimensional irreducible representations of
  $\cU_q(gl(2))$}

\medskip
All the finite dimensional simple modules over $\cU_q(gl(2))$ are of
course cyclic. 
We call \typeA\  representations those that are
deformation of classical representations, and \typeB\ the
others. 
Knowing the ${\cal U}_{q}(sl(2))$ case, 
we only need to add a parameter related with 
the value of the $U(1)$ generator $k_1 k_2^2$. 
This parameter may be provided as $\lambda_1\lambda_2^2$ 
(value of $k_1 k_2^2$) and a sign,
or simply by the value $\lambda_2$ of $k_2$ on a given vector. 
The finite dimensional irreducible 
representations of ${\cal U}_{q}(sl(2))$ are \cite{RA}:

\begin{itemize}
\item \typeA: 
  Usual (nilpotent) representations, where
  $k_1^{2l}=1$, $e_1^{l'}=0$, $f_1^{l'}=0$, characterized by their dimension
  $N=1,\cdots,l'$ and a sign $\omega$. The highest weight $\lambda_1$ 
  is $\lambda_1=q^{N-1}$.  
  These representations are given explicitly
  in Appendi<x A (\ref{eq:gl2_nil}). 
  The representation of dimension $l'$ plays a special role. It is in fact
  in the intersection of this case and the following. 
\item \typeB:
  Coloured (nilpotent) representations, with still $e_1^{l'}=0$,
  $f_1^{l'}=0$, characterized by their
  highest weight $\lambda_1$, a continuous parameter. Their dimension is
  $l'$. They are also described by (\ref{eq:gl2_nil}).
\item \typeB:
  Periodic and semi-periodic representations, explicitly
  given in (\ref{eq:gl2_per}). These representations have dimension
  $l$. 
  They depend on four complex parameters corresponding to the values
  of the three central elements $k_1^{l'}$, 
  $e_1^l$, $f_1^l$ and one discrete
  parameter corresponding to the value of the quadratic Casimir
  $\cC_{\cU_q(gl(2))}$ of ${\cal U}_{q}(gl(2))$ and related
  to the former through the relation (\ref{eq:rel_gl2}). 
  The ${\cal U}_{q}(gl(2))$-representation 
  is also completely characterized by the parameters
  $y=\varphi^l$, $\beta$, $\lambda_1$ and $\lambda_2$  appearing in 
  \begin{equation}
    \begin{array}{ll}
      f_{1}^l v_0 = \varphi^l v_0 \qquad &
      f_1 e_1 v_0 = \beta v_0 \\
      k_1 v_0 = \lambda_1 v_0 = q^{\mu_1}v_0
      \qquad  &
      k_2 v_0 = \lambda_2 v_0 = q^{\mu_2}v_0
    \end{array}
    \label{eq:gl2_per_short}
  \end{equation}
\end{itemize}

The existence of periodic irreducible representations has the
following consequence: the primitive ideals defined as the kernels of
these representations are not the annihilator of the irreducible
quotient of some Verma module, unlike in the case of classical
(super)algebras \cite{Duflo, MussonPrim}.

\subsection{Complete sets of representations of $\cU_q(sl(2))$}

We prove that a set of generic (periodic) representations
corresponding to an open subset of the set of parameters builds a
complete set, in the following sense: if an element of $\cU_q(sl(2))$
acts as $0$ on all the representations of this set, then it is the $0$
element of $\cU_q(sl(2))$. This terminology was used in \cite{ABP},
where the authors found complete sets of finite dimensional 
irreducible representations
of the classical $sl(2)$ and $sl(2|1)$. For quantum groups at roots of
unity, we shall obtain rather different results. 

Let $\cR\in \cU_q(sl(2))$ be such that it vanishes on a set $\Omega$
of representations. 
Let $q^{2(t-r)}$ be the $q$-grading  
of an element $f_1^r k_1^s e_1^t$. We have 
$k_1 \left( f_1^r k_1^s e_1^t \right) = q^{2(t-r)} \left( f_1^r k_1^s
e_1^t \right) k_1$. Any element of $\cU_q(sl(2))$ is a sum of terms of
given grading since the $f_1^r k_1^s e_1^t$ form a basis of
$\cU_q(sl(2))$. 
We write $\cR=\sum_{d=0}^{l'-1}\cR_d$, where the grading of $\cR_d$ is
$q^{2d}$. 
Commuting $\cR$ with $k_1$ shows that all the $\cR_d$ vanish
separately on the representations in $\Omega$. The same is true for
each $f_1^d\cR_d$. Since $\cU_q(sl(2))$ contains no zero divisor
\cite{KasselBook}, the vanishing of $f_1^d\cR_d$ in $\cU_q(sl(2))$ is
equivalent to that of $\cR_d$. 
Hence, to prove
that $\Omega$ is a complete set of representations, we only have to
show that the only element of $\cU_q(sl(2))$ commuting with $k_1$ and
acting as $0$ on all representations of $\Omega$ is $0$.

Let $\cR$ be an element of $\cU_q(sl(2))$ with grading $1$, 
and $n,n'\in\NN$
such that $f_1^n k_1^{n'}\cR = \sum_i a_i f_1^{r_i} k_1^{s_i}
e_1^{t_i}$ has 
only terms with $r_i - t_i\in l'\NN$ and $s_i\in \NN$. Then  $f_1^n
k_1^{n'}\cR$ can be written 
as a polynomial in $f_1^{l'}$,
$k_1$,  and $f_1e_1$, which commute among themselves. 
The value of this polynomial on
the vector $v_0$ of the representation (\ref{eq:gl2_per})
is the same polynomial evaluated on the scalars 
$\varphi^{l'}$, $\lambda_1$ and $\beta$. 
If $\Omega$ is a set of representations corresponding to an open
subset of $\CC^3$ for the values of $\varphi^{l'}$, $\lambda_1$ and
$\beta$, and if $\cR$ vanishes on all the representations of $\Omega$,
then the polynomial vanishes identically in $\cU_q(sl(2))$, and hence
$\cR=0$ as an element of  $\cU_q(sl(2))$. 
We then have the following

\begin{prop}
  A set of generic (periodic) representations
  corresponding to an open subset of the set of values for the 
  parameters is a
  complete set of representations. 
  \label{prop:complete_sl2}
\end{prop}

{\em Remark 1:} an element of $\cU_q(sl(2))$ that vanishes on all
\typeA\ modules, or even on all nilpotent or semiperiodic 
modules, is not necessarily
$0$ in $\cU_q(sl(2))$ (take simply $\cR = e_1^l f_1^l$). So a complete
set of representations should include periodic ones. 

{\em Remark 2:} suitably choosen 
infinite sets of periodic representations
(not necessarily corresponding to an open set of values of the
parameters) can also be complete. 

\section{Classification of finite dimensional 
  irreducible representations of ${\cal
    U}(sl(2|1))$ \label{sect:classification}}

Let $V_0$ an $N$-dimensional 
irreducible $\uqggo$-module, that we extend to a
$\uqggp$-module by the requirement that $e_2 V_0 = 0$.

Let $M'$ be the induced module $\cU_q(sl(2|1)) \otimes_{\uqggp}
V_0$. Then 
\begin{equation}
M' = V_0 \oplus f_2 V_0 \oplus f_3 V_0 \oplus f_2 f_3 V_0 \;.
\label{eq:M}
\end{equation}
The subspaces $(f_2 V_0 \oplus f_3 V_0)$ and $f_2 f_3 V_0$ are 
representations of ${\cal U}_{q}(\ggo)$ with the same value for 
central elements $k_1^{l'}$, $k_2^l$, 
$e_1^l$, $f_1^l$ as for $V_0$. If we write the value of quadratic Casimir
$\cC_{\cU_q(gl(2))}$ of $\cU_q(gl(2))$ as $\xi+\xi^{-1}$, then 
its eigenvalues on the different subspaces are
\begin{equation}
\begin{array}{lll}
\Bigg.\mbox{Subspace}\qquad\qquad\qquad & \cC_{\cU_q(gl(2))} & \\[.3cm]

V_0~: & \xi+\xi^{-1} & \\
(f_2 V_0 \oplus f_3 V_0)~: & 
q\xi+q^{- 1}\xi^{-1} \;,\quad& 
q^{- 1}\xi+q\xi^{-1}
\\
f_2 f_3 V_0~: & \xi+\xi^{-1} &
\end{array}
\end{equation}

\medskip

The elements $f_{2}^{\rho}f_{3}^{\sigma}f_{1}^{p}$, for $p\in \NN$,
$\rho\in\{0,1\}$ and $\sigma\in\{0,1\}$ build a
Poincar{\'e}--Birkhoff--Witt basis of the subalgebra 
$\cU^-$ generated by $f_1$ and
$f_2$. 
The elements $e_{1}^{p'}e_{3}^{\sigma'}e_{2}^{\rho'}$, for $p'\in \NN$,
$\rho'\in\{0,1\}$ and $\sigma'\in\{0,1\}$ build a
Poincar{\'e}--Birkhoff--Witt basis of the subalgebra 
$\cU^+$ generated by $e_1$ and
$e_2$. Together with the basis 
$k_1^{s_1} k_2^{s_2}$ (with
$s_i\in\ZZ$) for the Cartan subalgebra, this provide a basis for 
$\cU_q(sl(2|1))$.

Let $w_{0,0,0}$, $w_{0,0,1},\cdots, w_{0,0,N-1}$, be a basis of
$V_0$. Then it follows from the definition of $V_0$ and of the 
Poincar{\'e}--Birkhoff--Witt basis of $\cU_q(sl(2|1)$ given above 
that the vectors 
$f_{2}^{\rho}f_{3}^{\sigma} w_{0,0,p}$
$\rho, \sigma \in\{0,1\}$, $p\in\{0,N-1\}$ 
build a basis of $M'$. 
In particular 
\begin{equation}
  \dim M' = 4 N\;,
  \label{eq:dimM'}
\end{equation}
i.e. four times the dimension of $V_0$.

Since the dimension $N$ of $V_0$ is bounded by $l$, we already know
that the dimension of a simple $\cU_q(sl(2|1))$-module is bounded by
$4l$. Since nilpotent representations of $\uqggo$ have
dimension less or equal to $l'$, 
the dimension of nilpotent representations of
$\cU_q(sl(2|1))$ is bounded by $4l'$.

\subsection{Usual (\typeA) representations}

We now start from a $\uqggo$-module $V_0$ which is the $q$-deformation
of a classical module. Let $N$ be its dimension ($1\le N\le l'$).

The module $M'$ is then a highest weight module with highest weight
vector $w_{0,0,0}$ on which 
\begin{equation}
\begin{array}{ll}
  e_{1}\;w_{0,0,0}=0\;, &  e_{2}\;w_{0,0,0}=0 \;, \\
  k_{1}\;w_{0,0,0}=\lambda_1 w_{0,0,0}\;, \qquad
  & k_{2}\;w_{0,0,0}=\lambda_2 w_{0,0,0} \;
  \label{eq:hwv}
\end{array}
\end{equation}
with $\lambda_1 = \omega q^{N-1}$, $\omega=\pm 1$.

The Casimir operators $\cC_p$ have the following scalar value on $M'$: 
\begin{equation}
  {\cal C}_{p}= (q-q^{-1})^2 \lambda_{1}^{2p-1}\lambda_{2}^{4p-2}
  [\mu_{2}][\mu_{1}+ \mu_{2}+1] 
  \label{eq:Cp_valeur_nil}
\end{equation}
where, again, $q^{\mu_i}\equiv \lambda_i$.

A basis of $M'$ is given by 
\begin{equation}
  w_{\rho,\sigma,p}=f_{2}^{\rho}f_{3}^{\sigma}f_1^p w_{0,0,0}\;, 
  \qquad\mbox{with}\qquad
  \cases{ 
    \rho,\sigma\in\{0,1\}\cr
    p\in\{0,\cdots,N-1\}\;. \cr}
  \label{eq:base_M'1}
\end{equation}
By convention, we set 
\begin{equation}
  w_{\rho,\sigma,N} \equiv 0 \;.
\end{equation}

A non zero 
vector in a representation is called singular if it is annihilated
by $e_1$ and $e_2$ and is contained in a proper subrepresentation. 
Any submodule of $M'$ contains a singular vector for $M'$.
Indeed, any submodule of $M'$ has its own
$\uqggo$-submodule annihilated by $e_2$. This last module is also of
\typeA\ because this property is determined by the scalar value of the
central elements, which are determined by $V_0$.
The module $M'$ is simple if, and only if it contains no singular
vector $v_s\neq 0$.

\begin{lemma}
  The non-vanishing of the Casimir operators $\cC_p$ is a sufficient 
  condition for $M'$ to be simple.
\end{lemma}
The comparison of the values of the Casimir operators on the highest
weight vector and on the singular vector indeed shows that 
\begin{equation}
  [\mu_2][N + \mu_2] = 0
  \label{eq:Cp=0}
\end{equation}
is a necessary condition for the existence of a singular vector (which
cannot be in $V_0$ since $V_0$ is a simple $\uqggo$-module). This
condition amounts to the vanishing of all the Casimir operators
$\cC_p$. 
We shall see that this is actually a necessary and sufficient
condition for the simplicity of $M'$. 

\subsubsection{Typical \typeA\ representations}

\begin{prop}
If (\ref{eq:Cp=0}) is not satisfied, the module $M'$ is simple. It has
dimension $4N$.  
Its explicit expression is given in (\ref{eq:repr_nil}). It is called 
typical. 
\end{prop}
{\bf Proof:} If (\ref{eq:Cp=0}) is not satisfied, $M'$ contains no
singular vector. 

For $N=1,\cdots,l'-1$, the subspace $f_2V_0 \oplus f_3V_0$ is the
direct sum of $\uqggo$-modules characterized by the dimensions $N\pm
1$ and sign $\omega$. For $N=l'$, $f_2V_0 \oplus f_3V_0$ is an
indecomposable  $\uqggo$-module which is  
isomorphic to the tensor product of
$V_0$ with the spin $1/2$ representation, and which contains the
$\dim=l'-1$ (sign$=\omega$) simple sub-$\uqggo$-module.

\subsubsection{Atypical \typeA\ representations}
We now consider the case $[\mu_2][N+\mu_2]=0$
(i.e. $(\lambda_2^2-1)(\lambda_2^2 - q^{-2N})=0$). We will prove the
following:
\begin{prop}
  If the Casimir operators $\cC_p$ vanish on $M'$, there exists a maximal
  submodule $M''$ of $M'$. The quotient space $M=M'/M''$ is a simple
  module, called atypical. We can consider three cases:
  \begin{itemize}
  \item If $[\mu_2]=0$ and $[N+\mu_2]\neq 0$, then $\dim M=2N-1$.
  \item If $[\mu_2] \neq 0 $ and $[N+\mu_2]=0$ then $\dim M = 2N+1$.
  \item If $[\mu_2]= 0 $ and $[N+\mu_2]=0$ (and hence $N=l'$) then
    $\dim M=2l'-1$.
  \end{itemize}
\end{prop}

{\bf Proof:}

\paragraph{Atypical \typeA\ representations with $[\mu_2]=0$ and
  $[N+\mu_2]\neq 0$. }

In this case, the vector $f_2 w_{0,0,0}=w_{1,0,0}$ is a singular
vector. The action of $\cU_q(sl(2|1))$ on it generates a
$2N+1$-dimensional submodule $M''$ spanned by 
\begin{equation}
\begin{array}{ll}
  f_1^p w_{1,0,0} = q^{-p} w_{1,0,p} -q^{-1}[p]w_{0,1,p-1} 
  \;, \qquad 
  & p=0,\cdots ,N \\
  f_1^p f_3 w_{1,0,0} = -q^{-1} w_{1,1,p}
  \;, 
  & p=0,\cdots ,N-1  
  \label{eq:M''1}
\end{array}
\end{equation}
This submodule is maximal.
Quotienting $M'$ by $M''$ provides a $2N-1$-dimensional simple module
$M$, the expression of which is given in (\ref{eq:repr_nil_atyp1}). 

\paragraph{Atypical \typeA\ representations with $[\mu_2]\neq0$ and
  $[N+\mu_2]= 0$ }

Looking by direct computation for a singular vector, we see that
$N=1$, $[1+\mu_2]=0$
is a particular case: it is the only case of existence of a singular
vector in $f_2 f_3 V_0$ (one-dimensional in this case). The singular
vector is $w_{1,1,0}= f_2 f_3 w_{0,0,0}$. It generates only 
$f_2 f_3 V_0$ as $\cU_q(sl(2|1))$-submodule. The quotient 
$
M' / f_2 f_3 V_0
$
is three dimensional. It is actually the $q$-deformed three
dimensional atypical fundamental representation. 

If $N\in\{2,\cdots,l'-1\}$, there is a singular vector given by 
\begin{equation}
  v_{s} = \lambda_1 q w_{1,0,1} + [\mu_1] w_{0,1,0} \;.
  \label{eq:vec_sing}
\end{equation}
It generates the $2N-1$-dimensional 
maximal submodule $M''$ spanned by 
\begin{equation}
\begin{array}{ll}
  f_1^p v_s =  \lambda_1 q^{1-p} w_{1,0,p+1} + [\mu_1-p] w_{0,1,p}
  \;, \qquad
  & p=0,\cdots ,N-2 \\
  f_1^p f_2 v_s =  [\mu_1] w_{1,1,p}
  \;, 
  & p=0,\cdots ,N-1 
  \label{eq:M''2}
\end{array}
\end{equation}
The quotient $M=M'/M''$ is a $2N+1$-dimensional simple module,
explicitly given in (\ref{eq:repr_nil_atyp2})

\paragraph{Atypical \typeA\ representations with $N=l'$ }
If $[\mu_2]=0$ and $N=l'$, the vector $w_{1,0,0}$ is singular. 
The submodule it generates is similar to (\ref{eq:M''1}), except that
now 
$f_1^{l'} w_{1,0,0}=0$. However, the vector $w_{0,1,l'-1}$ is
subsingular, i.e. its image by $e_1$ and $e_2$ is contained is the
submodule generated by $w_{1,0,0}$. It belongs to the maximal
submodule $M''$ of $M'$. 
Note that $f_1 w_{1,0,0}\in M''$ is also singular.  
The submodule $M''$ has dimension $2l'+1$ and
$M=M'/M''$ has dimension $2l'-1$. It is also described by
(\ref{eq:repr_nil_atyp1}). 

\qed 

\subsection{Nilpotent \typeB\ representations} 
We now consider the case where $V_0$ is a \typeB\ nilpotent 
$\uqggo$-module, of dimension $N=l'$, with two parameters $\lambda_1$
and $\lambda_2$. We assume $[\mu_1+1]\neq 0$ since this case was
treated as \typeA. 
As in the \typeA\ case we consider the induced module $M'$, on which
(\ref{eq:hwv}) applies. A basis for $M'$ is also given by
(\ref{eq:base_M'1}) with $N=l'$. 
We also have
\begin{prop}
  Nilpotent \typeB\ representations fall into two classes
  \begin{itemize}
  \item If $[\mu_2][\mu_1 + \mu_2 + 1] \neq 0$, i.e. $\cC_p\neq 0$,
    then $M'$ is simple. Its dimension is $4l'$ and the parameters are
    $\lambda_1$ and $\lambda_2$. Its explicit expression is given in
    (\ref{eq:repr_nil})
    (typical case). 
  \item If $[\mu_2][\mu_1 + \mu_2 + 1] = 0$, i.e. $\cC_p= 0$, then
    $M'$ has a maximal submodule $M''$ of dimension $2l'$. Then
    $M=M'/M''$ has dimension $2l'$ (atypical case). 
  \end{itemize}
\end{prop}
{\bf Proof:} As in the \typeA\ case, there is no singular vector if
the $\cC_p$ do not vanish. 
Suppose now that $[\mu_2][\mu_1 + \mu_2 + 1] = 0$ 
We can separate
this case into two subcases, according to which term of the product
vanishes (both terms cannot vanish simultaneously, since $[\mu_1-p+1]\neq
0$ for any integer $p$ in \typeB\ $\uqggo$-modules). 
\begin{itemize}
\item If $[\mu_2]=0$, the vector $w_{1,0,0}$ is singular. It generates
  the submodule $M''$ given as in (\ref{eq:M''1}) with $N=l'$, except
  that now
  $f_1^{l'} w_{1,0,0}=0$. Then $\dim M''=2l'$. The quotient module
  hence has dimension $2l'$. It is described by
  (\ref{eq:repr_nil_atyp1}).  
\item If $[\mu_1 + \mu_2 + 1] = 0$, then there is a singular vector
  given by (\ref{eq:vec_sing}). It generates the submodule $M''$ given
  as in (\ref{eq:M''2}) with $N=l'$, except
  that now $f_1^{l'-1}v_s\neq 0$ also belongs to $M''$, so that $\dim
  M''=2l'$. Again, $\dim M=2l'$ and $M$ is described by
  (\ref{eq:repr_nil_atyp2}).  
\end{itemize}

\qed

\subsection{Periodic and semi-periodic \typeB\ representations} 

Let us now consider the case when $V_0$ is a periodic or
semi-periodic $\uqggo$-module, i.e. with non vanishing (scalar)
value of the central element $f_{1}^{l}$, 
\begin{equation}
  f_{1}^l = \varphi^l  \id \;.
  \label{eq:periodic}
\end{equation}
In $\cU_q(sl(2|1))$, 
$f_{1}^{l}$ is central too, so (\ref{eq:periodic}) also holds in $M'$.

The value of the central element $e_1^l$ will be a free parameter
(possibly zero for semi-periodic representations).
One would get the representations with a vanishing value for $f_1^l$ and a
non vanishing 
value for $e_1^l$, using the automorphism of $\cU_q(sl(2|1))$ given by 
\begin{equation}
  \begin{array}{ll}
    \psi(e_i) = f_i \qquad & \psi(f_i) = e_i \\
    \psi(k_1) = k_1^{-1}  \qquad & \psi(k_2) = -k_2^{-1}
  \end{array}
  \label{eq:automophism}
\end{equation}

The module $M'$ is actually characterized by the following actions on
a vector $w_{0,0,0}$ of $V_0$:
\begin{equation}
  \begin{array}{ll}
    f_{1}^l w_{0,0,0} = \varphi^l w_{0,0,0} \qquad &
    f_1 e_1 w_{0,0,0} = \beta w_{0,0,0} \\
    k_1 w_{0,0,0} = \lambda_1 w_{0,0,0} = q^{\mu_1}w_{0,0,0} \qquad  &
    k_2 w_{0,0,0} = \lambda_2 w_{0,0,0} = q^{\mu_2}w_{0,0,0} \\
  \end{array}
  \label{eq:per_short}
\end{equation}
Those values determine the values of $e_1^l$ (using (\ref{eq:rel_gl2}))
and of $\cC_p$
\begin{equation}
  {\cal C}_{p}= (q-q^{-1})^2 \lambda_{1}^{2p-1}\lambda_{2}^{4p-2}
  \left([\mu_{2}][\mu_{1}+ \mu_{2}+1] - \beta \right)
  \label{eq:Cp_valeur}
\end{equation}

A basis of $M'$ is given by 
\begin{equation}
  w_{\rho,\sigma,p} \equiv \varphi^{-\sigma-p}
  f_{2}^{\rho}f_{3}^{\sigma}f_1^p w_{0,0,p}\;, 
  \qquad\mbox{with}\qquad
  \cases{ 
    \rho,\sigma\in\{0,1\}\cr
    p\in\{0,\cdots,l-1\}\;. \cr}
  \label{eq:base_M'2}
\end{equation}

\begin{prop}
  For periodic and semi-periodic representations, the following
  alternative holds:
  \begin{itemize}
  \item 
    {\bf (i).} If 
    $
    [\mu_{2}][\mu_{1}+ \mu_{2}+1] - \beta \neq 0 \;, 
    $
    the module $M'$ is irreducible and its dimension is 
    equal to $4\;l$. It is described explicitly in equation 
    (\ref{eq:repr_per}).
  \item
    {\bf (ii).} If $[\mu_{2}][\mu_{1}+ \mu_{2}+1] - \beta = 0$, 
    the module $M'$ is not simple. 
    It has a submodule $M''$ of dimension $2l$ 
    and the factor space
    $  M' / M'' $
    is an irreducible module of dimension $2\;l$, explicitly given by
    equation 
    (\ref{eq:repr_per_aty}). 
  \end{itemize}  
\end{prop}

The cases (ii) corresponds to atypical periodic representations
and $[\mu_{2}][\mu_{1}+ \mu_{2}+1] = \beta$ is the condition for the
vanishing of the Casimir operators $\cC_p$ on $M'$.

{\bf Proof.} 
By direct computation, we check that $[\mu_{2}][\mu_{1}+ \mu_{2}+1] -
\beta = 0$ is the necessary and sufficient condition for the existence
of a vector (not belonging to $V_0$), annihilated by both $e_2$ and
$e_3$. This vector then belongs to $f_2 V_0 \oplus f_3 V_0$ and it
generates a $2l$-dimensional subspace 
spanned by the 
vectors $w_{1,1,p}$ and 
$[\mu_{2}+p+1]w_{0,1,p}-\lambda_2^{-1}q^{-p}w_{1,0,p+1}$ 
 for $p\in\{0,\cdots,l-1\}$.
The quotient of $M'$ by this submodule is simple.

\section{Complete sets of representations of $\cU_q(sl(2|1))$ 
  \label{sect:complete_set}}
\begin{prop}
  A set of typical periodic 
  representations corresponding to an open subset of
  the set of values of the parameters is a complete set of
  representations. 
\end{prop}
{\bf Proof:} Let $\Omega$ be a set of representations, and
$\cR\in\cU_q(sl(2|1))$ such that $\cR$ vanishes on all the
representations 
of $\Omega$. As for $\cU_q(sl(2))$, we can restrict ourselves to the case
where $k_i \cR k_i^{-1} = q^{d_i} \cR$ for given gradings $d_i$ ($i=1,2$).

We have in fact to consider five cases, according to the possible
gradings with respect to $k_1 k_2^2$. All the possible values for
$d_1+2d_2$ are actually $-2,-1,0,1,2$ (This is due to the fact that
the squares of fermionic generators vanish, and it can also be read
from the 
Poincar{\'e}--Birkhoff--Witt basis). 

\begin{equation}
  \begin{array}{ll}
    d_1+2d_2 = -2 \qquad& \cR^{(-2)} = \cR_1 e_3 e_2 \\
    d_1+2d_2 = -1 & \cR^{(-1)} = \cR_2 e_2 + \cR_3 e_3 + \cR_4 f_2 e_3 e_2
    + \cR_5 f_3 e_3 e_2 \\
    d_1+2d_2 =  0 & \cR^{(0)} = \cR_6 +\cR_7 f_2 e_2 +\cR_8 f_3 e_2 +\cR_9
    f_2 e_3 +\cR_{10} f_3 e_3 +\cR_{11} f_2 f_3 e_3 e_2 \\
    d_1+2d_2 =  1 & \cR^{(1)} = \cR_{12} f_2 + \cR_{13} f_3 + \cR_{14} f_2
    f_3 e_2 + \cR_{15} f_2 f_3 e_3\\
    d_1+2d_2 =  2 & \cR^{(2)} = \cR_{16} f_2 f_3 
  \end{array}
  \label{eq:cR_i}
\end{equation}
where the $\cR_i$ are elements of $\uqggo$. We have to prove that all
of them vanish. 
Since $\Omega$ is a set of representations corresponding to an open
subset of the set of values of the parameters, the representations of
$\uqggo$ given by the corresponding $V_0$ is a complete set. 
If we identify $V_0$ and $f_2f_3 V_0$ 
(as $\uqggo$-modules), we see that the vanishing 
of $\cR_1$ and $\cR_{16}$ results from this. 
Let us now consider $\cR^{(0)}$, the cases of $\cR^{(-1)}$ and
$\cR^{(1)}$ being simpler. Since $\cR^{(0)}e_3e_2 = \cR_6 e_3e_2$
act as zero on all the representations of $\Omega$, then $\cR_6=0$. 
Now, $\cR^{(0)}e_2 = (\cR_9 f_2 +\cR_{10} f_3 ) e_3 e_2$. 
This operator sends $f_2f_3 V_0$ into $f_2 V_0 \oplus f_3 V_0$ and is
supposed to act as zero.
Looking at
the explicit action of this operator on the vector $v_{1,1,p}$ and
using the fact that $f_2 V_0 \oplus f_3 V_0$ is generically a direct
sum of two inequivalent $\uqggo$-modules, we learn that
$\cR_9=\cR_{10}=0$. Multiplying $\cR^{(0)}$ on the right by $e_3$, we
then prove in a similar way that $\cR_7=\cR_8=0$. Finally, the proof
that $\cR_{11}=0$ mimics the proof of proposition
\ref{prop:complete_sl2}.

\section{Proof of the relation in the centre
  \label{sect:proof_rel}}

We now use a complete set of representation to prove the relation
\begin{eqnarray}
  \cP_l(\cC_1,\cdots,\cC_l) &\equiv &
  ({\cal C}_{1}+1)^{l}-1 + \sum_{m\geq 2 \;\; n\geq 0 \atop m+n \leq l}
  {\cal C}_{m}{\cal C}_{1}^{n}{l\over m-1}\biggl({m+n-1 \atop n+1}\biggr)
  \biggl({l-m \atop n}\biggr) \nonumber\\
  &=&
  \left( 1-z_1^2 z_2^2 \right)\left(z_2^2 -1\right) -
    (q-q^{-1})^{2l} 
    z_1^2 z_2^4 y_{1}x_{1} \;. 
  \label{eq:rel2_2}
\end{eqnarray}

On a typical \typeB\ periodic representation characterized by the
parameters $\lambda_1$, $\lambda_2$, $\varphi^l$ and $\beta$, the
value of $\cC_p$ is 
\begin{eqnarray}
  \cC_p &=& \left( \lambda_1 \lambda_2^{2} \right)^{2p-2} \cC_1
  \nonumber\\
         &=&  \lambda_1^{2p-1} \lambda_2^{4p-2}
  \left( (q\lambda_1\lambda_2-q^{-1}\lambda_1^{-1}\lambda_2^{-1})
    (\lambda_2-\lambda_2^{-1}) - (q-q^{-1})^2\beta 
  \right)
  \nonumber\\
  &=& \lambda_1^{2p-1} \lambda_2^{4p-2}
  \left( q\lambda_1\lambda_2^2 + q^{-1}\lambda_1^{-1}\lambda_2^{-2} 
    - (\xi+\xi^{-1}) 
  \right) \nonumber\\
  &=& \lambda_1^{2p-1} \lambda_2^{4p-2}
  \left(q^{1/2} \lambda_1^{1/2} \lambda_2 \xi^{1/2} 
    -  q^{-1/2} \lambda_1^{-1/2} \lambda_2^{-1} \xi^{-1/2} 
  \right)
  \nonumber\\
  &&\qquad\qquad \cdot
  \left(q^{1/2} \lambda_1^{1/2} \lambda_2 \xi^{-1/2} 
    -  q^{-1/2} \lambda_1^{-1/2} \lambda_2^{-1} \xi^{1/2} 
  \right)
  \label{eq:Cp_value}
\end{eqnarray}
where 
$
(q-q^{-1})^{-2}(\xi+\xi^{-1})\equiv (q-q^{-1})^{-2} (q\lambda_1 +
q^{-1}\lambda_1^{-1}) + \beta
$ 
is the value of the $\uqggo$
quadratic Casimir operator on the subspace $V_0$. 

The polynomial $\cP_l$ in (\ref{eq:rel2}) is such that, 
if we set
\begin{eqnarray}
  \cC_1 &=& \lambda_1 \lambda_2^{2} (x_1 - x_1^{-1})(x_2 - x_2^{-1}) 
  \nonumber\\
  \frac{x_2}{x_1} &=& \lambda_1 \lambda_2^{2}   
\end{eqnarray}
then
\begin{equation}
  \cP(\cC_1,\cdots,\cC_l) = \lambda_1^{l} \lambda_2^{2l}
  \left( x_1^{l} - x_1^{-l} \right)
  \left( x_2^{l} - x_2^{-l} \right)
\end{equation}
so that
\begin{equation}
  \cP(\cC_1,\cdots,\cC_l) = \lambda_1^{l} \lambda_2^{2l}
  \left( \lambda_1^l\lambda_2^{2l} + \lambda_1^{-l}\lambda_2^{-2l} 
    - (\xi^l+\xi^{-l}) \right) \;.
\end{equation}
Using the polynomial relation (\ref{eq:rel_gl2}) in $\uqggo$, we
identify $(\xi^l+\xi^{-l})$ with the value of 
$
(q-q^{-1})^{2l} f_1^{l} e_1^{l} + \left(k^{l} + k^{-l}\right)  
$
and we get the evaluation of the right hand side of (\ref{eq:rel2}) on
the representation. Since this is true for any typical periodic
representations, and since the set of those representations is
complete, the relation is true in the enveloping algebra. 

The existence of any other independent polynomial relation in the
centre would imply more relations among the parameters of the periodic
representations, so we also conclude that there is no other
independent relation. 

\appendix
\def\thesection{Appendix\ \Alph{section}}

\section{Finite dimensional irreducible representations of
  $\cU_q(gl(2))$}
\paragraph{Nilpotent modules of $\cU_q(gl(2))$}
\begin{equation}
  \begin{array}{ll}
  k_1 v_{p} = \lambda_1 q^{-2p} v_{p} \;, &\quad\mbox{for}\quad
  p\in\{0,\cdots,N-1\}\;, \\
  f_1 v_{p} = v_{p+1} & \quad\mbox{for}\quad
  p\in\{0,\cdots,N-2\}\;, \quad\mbox{and}\quad   f_1 v_{N-1} = 0 \;, \\ 
  e_1 v_{p} = [p][\mu_1 - p+1]   v_{p-1} \;, & \quad q^{\mu_1}\equiv
  \lambda_1 \\  
  k_2 v_{p} = \lambda_2 q^{p} v_{p} \;, &\quad\mbox{for}\quad
  p\in\{0,\cdots,N-1\}.   
  \end{array}
  \label{eq:gl2_nil}
\end{equation}

The dimension $N$ is the smallest non
negative integer satisfying $[N][\mu_1-N+1]=0$. 
For usual \typeA\ representations, $N\in\{1,\cdots,l'\}$ and 
the highest weight is related to $N$ by 
$\lambda_1 = \omega q^{N-1}$, with
$\omega = \pm 1$. 

For nilpotent \typeB\ representations 
$N=l'$ and $\lambda_1$ is a free parameter. 

If $N=l'$ and $\lambda_1=\pm q^{-1}$, the
representation is still the $q$-deformation of a classical one, but it
has $q$-dimension $[N]=0$. This case plays a special role. 

\paragraph{Periodic and semi-periodic modules of $\cU_q(gl(2))$}
\begin{eqnarray}
  k_1 v_{p} &=& \lambda_1 q^{-2p} v_{p} \;, \nonumber\\
  f_1 v_{p} &=& \varphi v_{p+1} \;, \nonumber\\
  e_1 v_{p} &=& \varphi^{-1} ([p][\mu_1 - p+1] + \beta)
  v_{p-1} \;, \nonumber\\ 
  k_2 v_{p} &=& \lambda_2 q^{p} v_{p} \;, 
  \label{eq:gl2_per}
\end{eqnarray}
with $p\in\{0,\cdots,l-1\}$, and $q^{\mu_i}\equiv\lambda_i$, without
defining $\mu_i$ itself.  These representations have no classical
analogue (\typeB).

\section{Finite dimensional irreducible representations of
  $\cU_q(sl(2|1))$}

The following relations are used to determine the action of the
generators on the representations:
\begin{eqnarray}
  &&f_{1}\;f_{2}^{\rho}f_{3}^{\sigma}f_{1}^{p}=
  q^{\sigma-\rho}f_{2}^{\rho}f_{3}^{\sigma}f_{1}^{p+1}-
  \rho (1-\sigma)q^{-\rho} f_{2}^{\rho-1}f_{3}^{\sigma+1}f_{1}^{p},
  \nonumber\\ 
  &&f_{2}\;f_{2}^{\rho}f_{3}^{\sigma}f_{1}^{p}=
  (1-\rho)\;f_{2}^{\rho+1}f_{3}^{\sigma}f_{1}^{p}, \\
  &&[e_{1},\;f_{2}^{\rho}f_{3}^{\sigma}f_{1}^{p}]=
  \sigma(1-\rho)(-1)^{\sigma}f_{2}^{\rho+1}f_{3}^{\sigma-1}f_{1}^{p}
  q^{h_{1}-2p+1}
  +[p]f_{2}^{\rho}f_{3}^{\sigma}f_{1}^{p-1}[h_{1}-p+1],
  \nonumber\\ 
  &&e_{2}\;f_{2}^{\rho}f_{3}^{\sigma}f_{1}^{p} 
    - (-1)^{\rho+\sigma}f_{2}^{\rho}f_{3}^{\sigma}f_{1}^{p} \; e_{2}=
  \rho \;f_{2}^{\rho-1}f_{3}^{\sigma}f_{1}^{p}[h_{2}+p+\sigma]+
  \sigma
  (-1)^{\rho}f_{2}^{\rho}f_{3}^{\sigma-1} f_{1}^{p+1}q^{-h_{2}-p} \;,
  \nonumber 
  \label{eq:ffff}
\end{eqnarray}
where $(p,\;\rho,\;\sigma)\in\NN\times\lbrace
0,\;1\rbrace\times\lbrace 0,\;1\rbrace$. 

\paragraph{Typical nilpotent modules}

\begin{eqnarray}
  &&k_{1}\;w_{\rho,\sigma,p}=\lambda_1
  q^{\rho-\sigma-2p}\;w_{\rho,\sigma,p}, 
  \nonumber\\ 
  &&k_{2}\;w_{\rho,\sigma,p}=\lambda_2
  q^{\sigma+p}\;w_{\rho,\sigma,p}, \nonumber\\ 
  &&f_{1}\;w_{\rho,\sigma,p}=
  q^{\sigma-\rho}w_{\rho,\sigma,p+1}-
  \rho (1-\sigma) q^{-\rho}w_{\rho-1,\sigma+1,p}, \nonumber\\
  &&f_{2}\;w_{\rho,\sigma,p}=
  (1-\rho)\;w_{\rho+1,\sigma,p}, \nonumber\\
  &&e_{1}\;w_{\rho,\sigma,p}= - 
  \sigma(1-\rho)\lambda_1 q^{-2p+1}w_{\rho+1,\sigma-1,p}
  +[p][\mu_{1}-p+1]w_{\rho,\sigma,p-1}, \nonumber\\
  &&e_{2}\;w_{\rho,\sigma,p}=
  \rho [\mu_{2}+p+\sigma]\;w_{\rho-1,\sigma,p}+
  \sigma (-1)^{\rho}\lambda_2^{-1}q^{-p}w_{\rho,\sigma-1,p+1}. 
  \label{eq:repr_nil}
\end{eqnarray}
with $(p,\;\rho,\;\sigma)\in \{0,\cdots,N-1\}\times\lbrace
0,\;1\rbrace\times\lbrace 0,\;1\rbrace$ in the left hand side and, by
convention,  $w_{\rho,\sigma,N} = 0$ in the right hand side. 
For \typeA\ modules, $q^{\mu_1}\equiv \lambda_1 = \omega q^{N-1}$. 
For \typeB\ nilpotent modules, $N=l'$ and $q^{\mu_1}\equiv \lambda_1$
is free. 

\paragraph{Atypical nilpotent modules; case $[\mu_2]=0$}

\begin{eqnarray}
  &&k_{1}\;w_{\sigma,p}=\lambda_1 q^{-\sigma-2p}\;w_{\sigma,p},
  \nonumber\\ 
  &&k_{2}\;w_{\sigma,p}=\varepsilon q^{\sigma+p}\;w_{\sigma,p}, \nonumber\\
  &&f_{1}\;w_{\sigma,p}=q^{\sigma}w_{\sigma,p+1}, \nonumber\\
  &&f_{2}\;w_{\sigma,p}=(1-\sigma)q^{p-1}[p]\;w_{\sigma+1,p-1}, \nonumber\\
  &&e_{1}\;w_{\sigma,p}=q^{-\sigma}[p][\mu_1+1-p-\sigma]w_{\sigma,
    p-1}, \nonumber\\ 
  &&e_{2}\;w_{\sigma,p}=
  \sigma \varepsilon q^{-p}w_{\sigma-1,p+1} \;.
  \label{eq:repr_nil_atyp1}
\end{eqnarray}
where $\sigma\;\in \lbrace 0,\;1\rbrace$. 
For \typeA\ representations, $p\in \{0,\cdots, N-1-\sigma\}$
and the dimension is $2\;N-1$. For \typeB\ representations, 
$p\in \{0,\cdots, l'-1 \}$
and the dimension is $2\;l'$.

\paragraph{Atypical nilpotent modules; case $[\mu_1+\mu_2+1]=0$}

\begin{eqnarray}
  &&k_{1}\;w_{\sigma,p}=\lambda_1 q^{-\sigma-2p}\;w_{\sigma,p},
  \nonumber\\ 
  &&k_{2}\;w_{\sigma,p}=\varepsilon \lambda_1^{-1}
  q^{\sigma+p-1}\;w_{\sigma,p}, \nonumber\\ 
  &&f_{1}\;w_{\sigma,p}=q^{\sigma}w_{\sigma,p+1}, \nonumber\\
  &&f_{2}\;w_{\sigma,p}=-(1-\sigma)\lambda_1^{-1}
  q^{p-2}[\mu_1-p+1]\;w_{\sigma+1,p-1}, \nonumber\\ 
  &&e_{1}\;w_{\sigma,p}=q^{-\sigma}[p+\sigma][\mu_1+1-p]w_{\sigma,
    p-1}, \nonumber\\ 
  &&e_{2}\;w_{\sigma,p}=
  \sigma \varepsilon \lambda_1 q^{-p+1}w_{\sigma-1,p+1} \;.
  \label{eq:repr_nil_atyp2}
\end{eqnarray}
where $\sigma\;\in \lbrace 0,\;1\rbrace$. 
For \typeA\ representations, $p\in \{-\sigma,\cdots, N-1\}$
and the dimension is $2\;N+1$. For \typeB\ representations, 
$p\in \{0,\cdots, l'-1 \}$
and the dimension is $2\;l'$.

\paragraph{Typical periodic modules}

The actions of the generators $e_{1}$, $e_{2}$, 
$f_{1}$ and $f_{2}$ on a typical periodic $M$ module 
are given by
\begin{eqnarray}
  &&k_{1}\;w_{\rho,\sigma,p}=\lambda_1
  q^{\rho-\sigma-2p}\;w_{\rho,\sigma,p}, 
  \nonumber\\ 
  &&k_{2}\;w_{\rho,\sigma,p}=\lambda_2
  q^{\sigma+p}\;w_{\rho,\sigma,p}, 
  \nonumber\\ 
  &&f_{1}\;w_{\rho,\sigma,p}= \varphi
  q^{\sigma-\rho}w_{\rho,\sigma,p+1}-\varphi
  \rho (1-\sigma) q^{-\rho}w_{\rho-1,\sigma+1,p}, \nonumber\\
  &&f_{2}\;w_{\rho,\sigma,p}=
  (1-\rho)\;w_{\rho+1,\sigma,p}, \nonumber\\
  &&e_{1}\;w_{\rho,\sigma,p}= - \varphi^{-1} 
  \sigma(1-\rho) \lambda_1 q^{-2p+1}w_{\rho+1,\sigma-1,p}
  +\varphi^{-1} \biggl([p][\mu_{1}-p+1]+
   \beta \biggr)w_{\rho,\sigma,p-1}, \nonumber\\
  &&e_{2}\;w_{\rho,\sigma,p}=
  \rho [\mu_{2}+p+\sigma]\;w_{\rho-1,\sigma,p}+
  \sigma (-1)^{\rho} \lambda_2^{-1} q^{-p}w_{\rho,\sigma-1,p+1} \;,
  \label{eq:repr_per}
\end{eqnarray}
with $(\rho,\sigma)\in \lbrace 0,\;1\rbrace^2$ and $p\in\{0,\cdots,l-1\}$.

\paragraph{Atypical periodic modules}
\begin{eqnarray}
  &&k_{1}\;{w}_{\sigma,p}=\lambda_1 q^{-\sigma-2p}\;{\tilde
    w}_{\sigma,p}, \nonumber\\ 
  &&k_{2}\;{w}_{\sigma,p}=\lambda_2 q^{\sigma+p}\;{\tilde
    w}_{\sigma,p}, \nonumber\\ 
  &&f_{1}\;{w}_{\sigma,p}=\varphi 
  q^{\sigma} {w}_{\sigma,p+1}, \nonumber\\
  &&f_{2}\;{w}_{\sigma,p}= 
  (1-\sigma) \lambda_2 q^{p-1}\;[\mu_{2}+p]\;{\tilde
    w}_{\sigma+1,p-1}, \nonumber\\
  &&e_{1}\;{w}_{\sigma,p}= \varphi^{-1} 
  q^{-\sigma}[p+\mu_{2}][\mu_{1}+\mu_{2}-p+1-\sigma]
  {w}_{\sigma,p-1}, \nonumber\\
  &&e_{2}\;{w}_{\sigma,p}= \sigma 
  \lambda_2^{-1} q^{-p} {w}_{\sigma-1,p+1},  
  \label{eq:repr_per_aty}
\end{eqnarray}
with $\sigma\in \lbrace 0,\;1\rbrace$ and $p\in\{0,\cdots,l-1\}$.


\begin{thebibliography}{99}

\bibitem {ABP} D. Arnal, H. Ben Amor and G. Pinczon, 
  {\sl The structure of $sl(2,1)$-supersymmetry: irreducible
    representations and primitive ideals,} 
  Pacific Jour. of Math. {\bf 165} No 1 (1994) 17.

\bibitem {Asutroisq} D. Arnaudon, 
  {\sl Periodic and flat irreducible representations of
    $SU(3)_{q}$.}
  Commun. Math. Phys. {\bf 134} (1990) 523.
  
\bibitem {ACF} D. Arnaudon, C. Chryssomalakos and L. Frappat,
  {\sl Classical and Quantum $sl(1|2)$ Superalgebras, Casimir
    Operators and Quantum Chain Invariants,} 
  q-alg/9503021, 
  Journ. of Math. Phys. {\bf 36/10} (1995).

\bibitem {DK}  C. De Concini and V.G. Kac, 
  {\sl Representations of quantum groups at roots of 1,}
  Progress in Math. {\bf 92} (1990) 471 (Birkh{\"a}user);
  {\sl Representations of quantum groups at roots of 1: reduction to the
    exceptional case,} preprint RIMS 792 (1991).
  
\bibitem {DKP}  C. De Concini, V.G. Kac and C. Procesi, 
  {\sl Quantum coadjoint action,}
  J. A.M.S. {\bf 5} (1992) 151;
  {\sl Some remarkable degenerations of quantum groups,}
  Commun. Math. Phys. {\bf 157} (1993) 405.
  
\bibitem {DobStAndrews} V.K. Dobrev, in Proc. {\sl Int. Group Theory
    Conference,} St Andrews, 1989, Vol. 1, Campbell and Robertson
  (eds.), London 
  Math. Soc. Lect. Notes Series 159, Cambridge University Press, 1991.
  
\bibitem {Dri} V.G. Drinfeld, {\sl Quantum Groups,} Proc.
  Int. Congress of Mathematicians, Berkeley, California, Vol.
  {\bf 1}, Academic Press, New York (1986), 798. 
  
\bibitem {Duflo} M. Duflo, 
  {\sl Sur la classification des id{\'e}aux primitifs dans l'alg{\`e}bre
    enveloppante d'une alg{\`e}bre de Lie semi-simple,}
  Ann. of Math. (2) {\bf 105} (1977) 107.

\bibitem {GCX} Ge Mo-Lin, Sun Chang-Pu and Xue Kang,
  {\sl New $R$-matrices for the Yang--Baxter equation associated with
    the representations of the quantum superalgebra $U_q(osp(1,2))$
    with $q$ a root of unity,}
  Phys. Lett. {\bf A 163} (1992) 176.

\bibitem {Jim} M. Jimbo, {\sl $q$-difference analogue of 
    $\ \cU (\cG )$ and the Yang Baxter equation,} Lett. Math. Phys. {\bf 10}
  (1985)  63.
  
\bibitem {Kaci} V. G. Kac, Funct. Anal. Appl. {\bf 9}, 263 (1975).

\bibitem {Kacii} V. G. Kac, Commun. Math. Phys. {\bf 53}, 31 (1977).

\bibitem{KacSuperRep} V.G. Kac, 
  {\sl Representations of classical Lie
    superalgebras}, in Lecture Notes in Mathematics {\bf 676},
  Springer-Verlag, Berlin, Heidelberg, New York, 1978.
  
\bibitem {KasselBook} C. Kassel,
  {\sl Quantum groups,} Graduate Texts in Math. {\bf 155},
  Springer-Verlag New York (1995). 

\bibitem {KerDipl} T. Kerler,
  {\sl Darstellungen der Quantengruppen und Anwendungen,}
  Diplomarbeit, ETH--Zurich, August 1989.
  
\bibitem {Kobayashi} K.-I. Kobayashi, {\sl Representation theory of
    $osp(1|2)_q$,} Z. Phys. {\bf C 59} (1993) 155.

\bibitem {MussonPrim} I. Musson,
  {\sl A classification of primitive ideals in the enveloping algebra
    of a classical simple Lie superalgebra,}
  Adv. Math. {\bf 91} (1992) 252.

\bibitem {PalStoa} T.D. Palev and N.I. Stoilova, 
  J. Phys. A Math. Gen. {\bf 26}  
  (1993) 5867.

\bibitem {PalStob} T.D. Palev and N.I. Stoilova, 
  {Unitarizable representations of the deformed parabose
    superalgebra $U_q(osp(1/2))$ at roots of 1,}
  J. Phys. A28 (1995) 7275.

\bibitem {PSVdJ} T.D. Palev, N.I. Stoilova and J. Van der Jeugt, {\sl
    Finite-dimensional representations of the quantum superalgebra
    $U_q[gl(n|m)]$ and related $q$-identities,} Commun. Math. Phys. {\bf
    166} (1994) 367. 
  
\bibitem {PalTol} T.D. Palev and V.N. Tolstoy, 
  {\sl Finite dimensional
    irreducible representations of the quantum superalgebra
    $\cU_q[gl(n|1)]$,} Commun. Math. Phys. {\bf 141} (1991) 549.
  
\bibitem {RA}  P. Roche and  D. Arnaudon,
  {\sl Irreducible representations of the quantum 
    analogue of $SU(2)$,}
  Lett. Math. Phys. {\bf 17} (1989) 295.
  
\bibitem {ScheuSuper} M. Scheunert, 
  {\sl The Theory of Lie
    Superalgebras,} Lecture  
  Notes in Mathematics, vol. 716 (Springer, Berlin, 1979).   

\bibitem{ScheuEig} M. Scheunert, {\sl Eigenvalues of Casimir operators
    for the general linear, the special linear, and the
    orthosympleptic Lie superalgebras,} 
  J. Math. Phys. {\bf 24} (1983) 2681. 
  
\bibitem {Zhanga} R. B. Zhang, 
  {\sl Finite dimensional irreducible
    representations of the quantum supergroup $U_q(gl(m/n))$,}
  J. Math. Phys. {\bf 34} (3), (1993) 1236. 

\bibitem {Zhangb} R. B. Zhang 
  {\sl Three manifold invariants arising
    from $U_q(osp(1/2))$,} 
  Mod. Phys. Lett. {\bf A9} (1994) 1453. 

\bibitem {Zhangc} R. B. Zhang 
  {\sl A two-parameter quantization of $sl(2|1)$ and its
    finite-dimensional  representations,}
  J. Phys. A: Math. Gen. {\bf 27} (1994) 817.

\end{thebibliography}
\end{document}